\documentclass[12pt,reqno]{amsart}
\usepackage{amssymb,mathrsfs,amsmath,amsfonts,amsthm,graphicx,pdfpages,csquotes}
\usepackage{xcolor}
\usepackage{fouriernc}
\usepackage[foot]{amsaddr}

\newcommand{\RCat}{\mathsf{RCat}}
\newcommand{\RMat}{\mathsf{RMat}}
\newcommand{\define}[1]{{\bf \boldmath{#1}}}
\newcommand{\op}{\mathrm{op}}

\newcommand{\namedcat}[1]{\mathsf{#1}}

\newcommand{\Cat}{\namedcat{Cat}}
\newcommand{\CAT}{\mathsf{CAT}}

\newcommand{\Set}{\namedcat{Set}}

\newcommand{\Open}{\mathsf{Open}}

\newcommand{\Mat}{\mathsf{Mat}}

\newcommand{\colim}{\mathrm{colim}}
\makeatletter
\newcommand*{\relrelbarsep}{.386ex}
\newcommand*{\relrelbar}{%
  \mathrel{%
    \mathpalette\@relrelbar\relrelbarsep
  }%
}
\newcommand*{\@relrelbar}[2]{%
  \raise#2\hbox to 0pt{$\m@th#1\relbar$\hss}%
  \lower#2\hbox{$\m@th#1\relbar$}%
}
\providecommand*{\rightrightarrowsfill@}{%
  \arrowfill@\relrelbar\relrelbar\rightrightarrows
}
\providecommand*{\leftleftarrowsfill@}{%
  \arrowfill@\leftleftarrows\relrelbar\relrelbar
}
\providecommand*{\xrightrightarrows}[2][]{%
  \ext@arrow 0359\rightrightarrowsfill@{#1}{#2}%
}
\providecommand*{\xleftleftarrows}[2][]{%
  \ext@arrow 3095\leftleftarrowsfill@{#1}{#2}%
}
\makeatother

\usepackage{subfiles}
\usepackage{stackengine}
\usepackage{mathtools}
\usepackage[utf8]{inputenc}
\usepackage{tikz}
\usepackage[a4paper,top=3cm,bottom=3cm,inner=3cm,outer=3cm]{geometry}

\usepackage{comment}

\usepackage{graphicx}
\usepackage{adjustbox}
\usepackage[all,2cell]{xy}\UseAllTwocells\SilentMatrices
\definecolor{darkgreen}{rgb}{0,0.45,0}
\usepackage[colorlinks,citecolor=darkgreen,urlcolor=gray,final,hyperindex,linktoc=page,pagebackref]{hyperref}

\usepackage[capitalize]{cleveref}
\crefname{equation}{}{}
\crefname{item}{}{}
\usepackage{enumerate}

\newtheorem*{thm*}{Theorem}
\theoremstyle{remark}
\newtheorem*{rmk*}{Remark}
\newtheorem*{lem*}{Lemma}
\theoremstyle{definition}
\newtheorem*{defn*}{Definition}
\newtheorem*{cor*}{Corollary}
\theoremstyle{definition}
\newtheorem*{examples*}{Examples}
\newtheorem{prop*}{Proposition}

\theoremstyle{plain}
\newtheorem{thm}{Theorem}[section]
\theoremstyle{plain}
\newtheorem{prop}[thm]{Proposition}
\theoremstyle{remark}

\theoremstyle{plain}
\newtheorem{lem}[thm]{Lemma}
\theoremstyle{plain}
\newtheorem{cor}[thm]{Corollary}
\theoremstyle{definition}
\newtheorem{defn}[thm]{Definition}
\theoremstyle{definition}

\newcommand{\maps}{\colon}

\usepackage{tikz}
\usepackage{tikz-cd}
\usetikzlibrary{backgrounds,circuits,circuits.ee.IEC,shapes,fit,matrix}

\usepackage{verbatim}
\usetikzlibrary{arrows,shapes}

\tikzstyle{simple}=[-,line width=2.000]
\tikzstyle{arrow}=[-,postaction={decorate},decoration={markings,mark=at position .5 with {\arrow{>}}},line width=1.100]
\pgfdeclarelayer{edgelayer}
\pgfdeclarelayer{nodelayer}
\pgfsetlayers{edgelayer,nodelayer,main}

\tikzstyle{none}=[inner sep=0pt]

\definecolor{lblue}{rgb}{0,250,255}
\tikzstyle{species}=[circle,fill=yellow,draw=black,scale=1.15]
\tikzstyle{transition}=[rectangle,fill=lblue,draw=black,scale=1.15]
\tikzstyle{inarrow}=[->, >=stealth, shorten >=.03cm,line width=1.5]
\tikzstyle{empty}=[circle,fill=none, draw=none]
\tikzstyle{inputdot}=[circle,fill=black,draw=black, scale=.25]
\tikzstyle{inputarrow}=[->,draw=purple, shorten >=.05cm]
\tikzstyle{simple}=[-,draw=black,line width=1.000]
\tikzstyle{dot}=[circle,fill=black,draw=black, scale=.4]
\tikzstyle{inarrow}=[->, >=stealth, shorten >=.03cm,line width=1.5]

\definecolor{joecolor(x11)}{rgb}{0.0, 0.5, 0.5}

\definecolor{purple(x11)}{rgb}{0.8, 0, 0.8}

\title{The Open Algebraic Path Problem}

\author{Jade Master}

\address{Department of Mathematics, University of California, Riverside, 900 University Avenue, 92521, USA}

\email{jmast003@ucr.edu}

\begin{document}

\begin{abstract}
    The algebraic path problem provides a general setting for shortest path algorithms in optimization and computer science. This work extends the algebraic path problem to networks equipped with input and output boundaries. We show that the algebraic path problem is functorial as a mapping from a double category whose horizontal composition is gluing of open networks. We introduce functional open matrices, for which the functoriality of the algebraic path problem has a more practical expression.
\end{abstract}

\maketitle
The algebraic path problem is a generalization of the shortest path problem to probability, computing, matrix multiplication, and optimization \cite{tarjan1981unified,foote2015kleene}. Let $R$ be the rig of positive real numbers $([0,\infty], \mathrm{min}, +)$. A weighted graph is regarded as a matrix weighted in $R$, and the shortest paths of this graph are computed as the transitive closure of this matrix. The algebraic path problem allows $R$ to vary, and gets solutions to other problems of a similar flavor also as the transitive closure of an adjacency matrix. Many popular shortest path algorithms can be extended to compute these transitive closures in a more general setting \cite{hofner2012dijkstra} and the algebraic path problem can also be implemented generically using functional programming \cite{dolan2013fun}.

The algebraic path problem deals only with closed systems, i.e. systems which are isolated from their surroundings. On the other hand, open systems are equipped with input and output boundaries, from which they can be composed to form larger and more complicated networks. A research program intiated by Baez, Courser, and Fong aims to provide a theoretical foundation for open systems using cospan formalisms \cite{fong2016algebra,structured}. For a category of networks $C$, Baez and Courser defined a symmetric monoidal double category which provides a syntax for composition of open systems in $C$ \cite{structured}. In Section \ref{openmat}, we set $C$ equal to $\RMat$, the category of matrices weighted in a quantale $R$, to obtain a symmetric monoidal double category $\Open(\RMat)$. The essence of this double category is gluing. Open $R$-matrices are represented as cospans with feet given by $0$-matrices. Given two such open $R$-matrices, take their pushout 
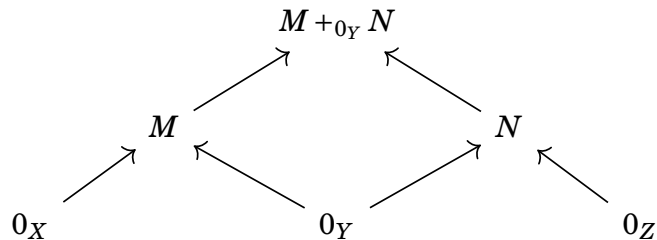
\[
\begin{tikzcd}
&   & M+_{0_Y} N &  & \\
& M \ar[ur] & & N\ar[ul] & \\
0_X\ar[ur] & & \ar[ul] 0_Y \ar[ur] & & \ar[ul] 0_Z
\end{tikzcd}
\]
to obtain an open $R$-matrix whose apex is synthesized from joining $M$ and $N$ along their shared boundary. This, along with the other data and structure of $\Open(\RMat)$, provide a syntax for manipulating open $R$-matrices. The axioms of a symmetric monoidal double category guarantee that this syntax is coherent. For example, the word problem for double categories is solvable in quadratic time \cite{wordproblem}.

$\RCat$, the category of $R$-enriched categories provide a choice of semantics for $R$-matrices, and can be expressed as $R$-matrices satisfying some regularity properties. In Section \ref{Kleene}, we show how the solution to the algebraic path problem forms the left adjoint $F$ of an adjunction
\[\begin{tikzcd}\RMat \ar[r,bend left,"F"] \ar[r,"\bot",phantom]& \RCat \ar[l,bend left,"U"] .\end{tikzcd} \]
which provides a mapping from the syntax of $R$-matrices to the semantics of $R$-categories. $R$-categories equipped with input and output boundaries form the horizontal morphisms of a symmetric monoidal double category $\Open(\RCat)$. In Section \ref{openalg}, we show how the algebraic path problem functor lifts to a symmetric monoidal double functor
\[\bigstar \maps \Open(\RMat) \to \Open(\RCat) \]
providing a coherent semantics for the syntax of \emph{open} $R$-matrices. This symmetric monoidal double functor explicates the way the solution to the algebraic path problem can built inductively from gluings of open $R$-matrices. The axioms of a symmetric monoidal double functor guarantee that this inductive process is coherent. 

This result is more theoretical than practical. However, there is a subclass of open $R$-matrices, functional open $R$-matrices, for which the theory provides useful insight. Functional open $R$-matrices are roughly open $R$-matrices where the inputs are all sources and the outputs are all sinks. In Section \ref{functional} we show that there is a strict double functor
\[\blacksquare \circ \bigstar_{fxn} \maps \Open(\RMat)_{fxn} \to \Mat_R \]\noindent  where $\Mat_R$ is a double category of $R$-matrices whose horizontal composition is matrix multiplication. This strict double functor gives a series of coherent compositional relationships for the algebraic path problem on functional open $R$-matrices based on matrix multiplication.

\section{The Algebraic Path Problem}\label{Kleene}

The algebraic path problem arises from the observation that various optimization problems can be framed in the same way by varying a sufficiently nice sort of rig. The level of generality for this work will be a commutative quantale, which is sufficient to guarantee existence and uniqueness of solutions to these optimization problems.
\begin{defn}
A \define{quantale} is a monoidal closed poset with all joins. Explicitly, a quantale is a a poset $R$ with a associative, unital, and monotone multiplication $\cdot \maps R \times R\to R$ such that
\begin{itemize}
    \item all joins, $\sum_{i \in I} x_i$, exist for arbitrary index set $I$ and,
    \item  $\cdot$ preserves all joins, i.e. 
    \[a \cdot \sum_{i \in I}  x_i = \sum_{i \in I} a \cdot x_i\]
    for all joins over an arbitrary index set $I$.
\end{itemize}
A quantale is commutative if its multiplication operation, $\cdot$, is commutative.
\end{defn}\noindent A motivating example of such a quantale is the poset $[0,\infty]$ with $+$ as its monoidal product and with join given by infimum. Note that this poset is equipped with the reverse of the usual ordering on $[0,\infty]$. Fong and Spivak show how the shortest path problem on this quantale computes the shortest paths between all pairs of vertices in a given $[0,\infty]$-weighted graph \cite[\S 2.5.3]{fong2019invitation}. Other motivating examples include the rig $([0,1],\mathrm{sup},\times)$ (whose algebraic path problem corresponds to most likely path in a Markov chain) and the powerset of the language generated by an alphabet (whose algebraic path problem corresponds to the language decided by a nondeterministic finite automata (NFA))\cite{foote2015kleene}.
\begin{defn}
For a commutative quantale $R$ and sets $X$ and $Y$, an \define{$R$-matrix} $M \maps X \to Y$ is a function $M \maps X \times Y \to R$. For $R$-matrices $M \maps X \to Y$ and $N \maps Y \to Z$, their matrix product $M N$ is defined by the rule
\[MN(i,k) = \sum_{j\in Y} M(i,j) N(j,k)\]
\end{defn}
\noindent If $R$ is a commutative quantale, $R$-matrices form a quantale as well. 
\begin{defn}
Let $\RMat(X)$ be the set of $X$ by $X$ matrices $M \maps X \times X \to R$. $\RMat(X)$ is equipped with the partial order $M \leq N$ if and only if $M(i,j) \leq N(i,j)$ for all $i,j \in X$.
\end{defn}
\begin{prop}
$\RMat (X)$ is a quantale with
\begin{itemize}
    \item join given by pointwise sum of matrices,
    \item and multiplication given by matrix product.
\end{itemize}
\end{prop}
\noindent The proof of this proposition is left to the reader. All the required properties of $\RMat(X)$ follow from the analogous properties in $R$.

A square matrix $M\maps X \times X \to R$ represents a complete $R$-weighted graph whose vertex set is given by $X$. 
\begin{defn} Let $M\maps X \times X \to R$ be a square matrix. A \define{vertex} of $M$ is an element $i\in X$. An \define{edge} of $M$ is a tuple of vertices $(a,b) \in X \times X$. A \define{path} in $M$ from $a_0$ to $a_n$ is a list of adjacent edges $p=\large((a_0,a_1),(a_1,a_2), \ldots,(a_{n-1},a_{n}), \large)$. The \define{weight} of $p$ is defined as the product
\[l(p) = \Pi_{i=0}^{n-1} M(a_i,a_{i+1}) \]
in $R$. For vertices $i,j \in X$, let
\[P_{ij}^M = \{\text{ paths in $M$ from $i$ to $j$ } \}  \]
\end{defn}Let $i$ and $j$ be vertices of a square matrix $M \maps X \times X \to R$. The algebraic path problem asks to compute the quantity
\[\sum_{p \in P_{ij}} l(p) \]
in the quantale $R$. If $R$ is the quantale $([0,\infty],\inf,+)$ then the weight of an edge $M_{ij}$ represents the distance between vertex $i$ and vertex $j$ and the weight of a path $l(p)$ represents the total distance traversed by $p$. Summing the weights of all paths between a pair of vertices corresponds to finding the path with the minimum weight. For example, the algebraic path problem asks to compute the length of the shortest path in the case when $R$ is $([0,\infty],\inf,+)$.

A more tractable framing of the algebraic path problem can be found by considering matrix powers. The entries of $M^2$ are given by 
\[M^2(i,j) = \sum_{l \in X} M(i,l)M(l,j)= \mathrm{inf}_{l \in X}\{ M(i,l)+ M(l,j)\}. \]
Because $M(i,l)$ and $M(l,j)$ represent the distance from $i$ to $l$ and from $l$ to $j$, this infimum computes the cheapest way to travel from $i$ to $j$ while stopping at some $l$ in between. More generally, the entries of $M^n$ for $n\geq 0$ represent the shortest paths between nodes of your graph that occur in exactly $n$ steps. To compute the shortest paths which can occur in any number of steps, we must take the infimum of the matrices $M^n$ over all $n \geq 0$. This pattern replicates for other choices of quantale. Therefore, the \define{algebraic path problem} seeks to compute 
\begin{equation}\label{pathproblem}
    F(M) = \sum_{n \geq 0} M^n
\end{equation}
where $M$ is an $R$-matrix. The following table summarizes some instances of the algebraic path problem for different choices of $R$. Fink provides an explanation of the algebraic path problems for $([0,\infty],\leq)$ and $\{T,F\}$ and Foote provides an explanation for the quantales $([0,1],\leq)$ and $(\mathcal{P}(\Sigma),\subseteq)$
\cite{fink1992survey,foote2015kleene}.
\smallskip
\begin{center}
\begin{tabular}{ |c|c|c|c| } 
\hline
\textbf{poset} & \textbf{join} & \textbf{multiplication} & \textbf{solution of path problem}\\
 \hline
 $([0,\infty],\geq)$ & $\inf$& $+$ & shortest paths in a weighted graph \\ 
 \hline
 $([0,\infty],\leq)$&$\sup$& $\inf$ & maximum capacity in the tunnel problem\\
 \hline
 $([0,1],\leq)$ & $\sup$ & $\times$ & most likely paths in a Markov process\\ 
 \hline
 $\{T,F\}$ & $\mathrm{OR}$ & $\mathrm{AND}$ & transitive closure of a directed graph \\ 
 \hline
 $(\mathcal{P}(\Sigma^*),\subseteq)$ & $\bigcup$ & concatenation & decidable language of a NFA \\
 \hline
\end{tabular}
\end{center}
\smallskip
Note that in this table, $\mathcal{P}(\Sigma^*)$ denotes the power set of the language generated by an alphabet $\Sigma$.

Equation (\ref{pathproblem}) is known to category theorists by a different name: the free monoid on $M$. Framing it in this way gives a categorical proof of existence and uniqueness of $F(M)$. A classic result from \cite[\S V11]{maclane} gives a construction of free monoids. MacLane's construction is defined as an adjunction into a category of internal monoids.
\begin{defn}
Let $(C, \otimes, I)$ be a monoidal category. A \define{monoid internal to $C$} is an object $A$ of $C$ equipped with morphisms
\[m \maps A \otimes A \to A \text{ and } i \maps I \to M \]
satisfying the axioms of associativity and unitality expressed as commutative diagrams. A \define{monoid homomorphism} from a monoid $A$ to a monoid $B$ is a morphism $f: A \to B$ in $C$ which commutes with the maps $m$ and $i$ of each monoid. Let \define{$\mathsf{Mon}(C)$} be the category where objects are monoids internal to $C$ and morphisms are their homomorphisms.
\end{defn}
\begin{prop}[MacLane]\label{Mac}Let $(C, \otimes, I)$ be a monoidal category with countable coproducts such that tensoring on both sides preserves these coproducts then there is an adjunction
\[\begin{tikzcd} \ \ \ C\ar[r,bend left,"F"] \ar[r,phantom,"\bot",pos=.6]& \ar[l,bend left,"U"] \mathsf{Mon}(C)\end{tikzcd} \]
whose left adjoint is given by the countable coproduct
\begin{equation}\label{freemon}
F(X) = \sum_{n \geq 0} X^n.
\end{equation}
\end{prop}
The poset $\RMat(X)$ when viewed as a category satisfies the hypotheses of Proposition \ref{Mac} and admits a free monoid construction.
\begin{prop}\label{fixedmon}
There is an adjoint pair
\[\begin{tikzcd}\RMat(X) \ar[r,bend left,"F_{X}"] &\mathsf{Mon}(\RMat(X)) \ar[l,bend left,"U_{X}"] \end{tikzcd} \]
where $F_{X}$ is the monotone map which produces the solution to the algebraic path problem on a matrix and $U_{X}$ is the natural forgetful map.
\end{prop}
\begin{proof}
Because $\RMat(X)$ is a quantale, it can be regarded as a monoidal category with all coproducts such that tensoring distributes over these coproducts. The result follows from applying Proposition \ref{Mac} and noticing that Equation \ref{freemon} matches Equation \ref{pathproblem} in the case when $C=\RMat$.
\end{proof}

\noindent Monoids internal to $\RMat(X)$ are $R$-enriched categories.
\begin{defn}
An $R$-category $C$ with object set $X$ consists of an element $C(x,y)$ in $R$ for every $x,y \in X$ such that
\begin{itemize}
    \item $1 \leq C(x,y)$ (the identity law),
    \item and $C(x,y) C(y,z) \leq C(x,z)$ (the composition law).
\end{itemize}
Let $\RCat(X)$ be the poset whose elements are $R$-enriched categories with object set $X$. For $R$-categories $C$ and $D$,
\[C \leq D \leftrightarrow C(i,j) \leq D(i,j) \quad \forall i,j \in X \]
\end{defn}
\begin{prop}
$\mathsf{Mon}(\RMat(X))$ is isomorphic to $\RCat(X)$, the poset of categories enriched in $R$ with object set $X$.
\end{prop}
\begin{proof}
The isomorphism in question assigns a matrix $M \maps X \times X \to R$ to the $R$-category with $\hom(x,y)=M(x,y)$. The identity law follows from the inequality $1 \leq M$ and the inequality $M^2 \leq M$ implies that for all $y \in X$,
\[\sum_{y \in X} M(x,y) M(y,z) \leq M(x,z) \] The composition law follows from the fact that any element of $R$ is less than a join which contains it.
\end{proof}

\noindent Proposition \ref{fixedmon} says that each matrix valued in $R$ has a unique, universally characterized solution to the algebraic path problem: namely the free $R$-category on that matrix. This adjunction can be extended to matrices over an arbitrary set. 
\begin{defn}\label{pushforward}
Let $f: X\to Y$ be a function and let $M: X \times X \to R$ be an $R$-matrix. Then the \define{pushforward} of $M$ along $f$ is the matrix $f_* (M) \maps Y \times Y \to R$ defined by 
\[f_*(M) (y,y') = \sum_{(x,x') \in (f \times f)^{-1}(y,y')} M(x,x'). \]
\end{defn}

\begin{defn}\label{matr}
Let $\RMat$ be the category where objects are square matrices $M \maps X \times X \to R$ on some set $X$ and where a morphism from $M \maps X \times X \to R$ to $N \maps Y \times Y \to R$ is a function $f \maps X \to Y$ satisfying 
\[f_*(M) \leq N. \]
Let $\RCat$ be the full subcategory of $\RMat$ consisting of matrices satisfying the axioms of an $R$-category.
\end{defn}

\begin{prop}
The free monoid construction of Proposition \ref{fixedmon} extends to an adjunction
\[\begin{tikzcd}\RMat \ar[r,bend left,"F"] & \RCat \ar[l,bend left,"U"] .\end{tikzcd} \]
\end{prop}

\begin{proof}
Let $A \maps \Set^{op} \to \Cat$ be the functor which sends a set $X$ to the poset $\RMat(X)$ regarded as a category and sends a function $f \maps X \to Y$ to the pushforward functor
\[f_* \maps \RMat(X) \to \RMat(Y).\]Analogously, let $B \maps \Set^{\op} \to \Cat$ be the functor which sends a set $X$ to the poset $\RCat(X)$ and sends a function $f$ to it's pushforward functor. The functors $F_{X}$ form the components of a natural transformation $\mathbf{F} \maps A \Rightarrow B$ and the functors $U_{X}$ form the components of a natural transformation $\mathbf{U} \maps B \Rightarrow A$. Furthermore, these natural transformations form an adjoint pair in the $2$-category $[\Set^{\op}, \Cat]$ of functors $\Set^{\op} \to \Cat$, natural transformations between them, and modifications. $\mathbf{F}$ and $\mathbf{U}$ are adjoint because an adjoint pair in $[\Set^{\op},\Cat]$ is a pair of natural transformations which are adjoint in each component. To summarize, we have a pair of adjoint natural transformations 
\[
\begin{tikzcd}
\Set^{\op}\, \ar[r,bend left=70,"A"{name=U}] \ar[r,bend right=70,"B"{name=D},swap] & \ar[Rightarrow,from=U, to=D,shorten <= 1.7ex, shorten >= 1.7ex,bend right,"\mathbf{F}"description] \ar[Rightarrow,from=D, to=U,shorten <= 1.7ex, shorten >= 1.7ex,bend right,"\mathbf{U}"description]\Cat
\end{tikzcd}
\]\noindent A restriction of the Grothendieck construction \cite{Borceux} is a 2-functor
\[\int \maps [\Set^{\op},\Cat] \to \CAT \]
where $\CAT$ is the 2-category of large categories. Because every 2-functor preserves adjunctions, the above diagram maps to an adjunction
\[ \begin{tikzcd}\int A \ar[r,bend left,"\int \mathbf{F}"] & \int B \ar[l,bend left,"\int \mathbf{U}"] .\end{tikzcd} \]
The result follows from the equivalences $\int A \cong \RMat$ and $\int B \cong \RCat$. The desired functors $F$ and $U$ are obtained by composing $\int \mathbf{F}$ and $\int \mathbf{U}$ with these equivalences.
\end{proof}

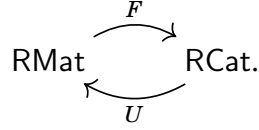
\begin{prop}
\[\begin{tikzcd}\RMat \ar[r,bend left,"F"] & \RCat \ar[l,bend left,"U"] .\end{tikzcd} \]
is an idempotent adjunction.
\end{prop}

\begin{proof}
Every adjunction between posets is idempotent. Therefore the smaller adjunctions $F_{X} \dashv U_{X}$ are idempotent. Because $F$ and $U$ are stitched together using these adjunctions, it is idempotent as well.
\end{proof}

\section{Open $R$-Matrices}\label{openmat}
$R$-matrices are made open by designating some of their vertices to be either inputs or outputs. In this section we show how these open $R$-matrices are composed by joining the output vertices of one to the input vertices of another and joining the data on the overlap. To define open $R$-matrices, we need a notion of a discrete weighted matrix on a set. The map sending a set to its discrete $R$-matrix is a functor and a left adjoint.
\begin{prop}\label{L}
Let $R \maps \RMat \to \Set$ be the functor which sends a weighted graph to its underlying set of vertices and sends a morphism to its underlying function. Then $R$ has a left adjoint
\[ 0 \maps \Set \to \RMat\]
which sends a set $X$ to the $R$-weighted graph
\[0_X \maps  X \times X \to Y \]
defined by $0_X (i,j) = 0$ for all $i$ and $j$ in $X$. $F$ sends a function $f \maps X \to Y$ to the morphism of $R$-matrices which has $f$ as its underlying function between vertices. 
\end{prop}
\begin{proof}
The natural isomorphism 
\[\phi\maps \RMat (0_X, G) \cong \Set (X, R(G))\]
is formed by noting that a morphism $ 0_X \to R(G)$ is uniquely determined by its underlying function on vertices and every such function obeys the inequality in Definition \ref{matr}.
\end{proof}
A weighted graph can be opened up to its environment by equipping it with inputs and outputs.
\begin{defn}Let $M: A \times A \to R$
An \define{open $R$-matrix} $M \maps X \to Y$ is a cospan in $\RMat$ of the form
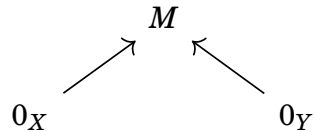
\[ \begin{tikzcd} & M & \\
 0_X \ar[ur] & & 0_Y \ar[ul]\end{tikzcd}\]
\end{defn}
\noindent The idea is that the maps of this cospan point to input and output nodes of the matrix $M$. Let  $M \maps X \to Y$ and $N \maps Y \to Z$
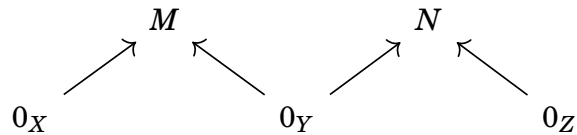
\[
\begin{tikzcd}
& M & & N & \\
0_X \ar[ur] & & 0_Y \ar[ul] \ar[ur] & & \ar[ul] 0_Z 
\end{tikzcd}
\]
be open $R$-matrices. The underlying sets of $M$ and $N$ form a diagram 
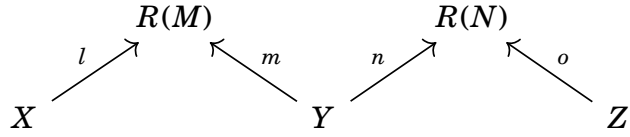
\[ \begin{tikzcd}
& R(M) & & R(N) & \\
X \ar[ur,"l"] & & Y \ar[ul,"m",swap] \ar[ur,"n"] & & \ar[ul,"o",swap] Z 
\end{tikzcd}\]
which generate a pushout
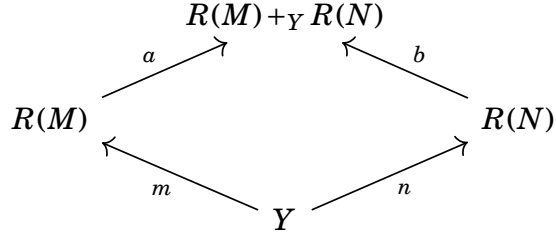
\[ 
\begin{tikzcd}
& R(M) +_Y R(N) & \\
R(M) \ar[ur,"a"] & & R(N) \ar[ul,"b",swap] \\
& Y \ar[ul,"m"] \ar[ur,"n",swap] & 
\end{tikzcd}
\]
The functions $a$ and $b$ of this pushout allow the matrices $M$ and $N$ to be compared on equal footing: the pushforwards $a_*(M)$ and $b_*(N)$ both have $R(M) +_Y R(N)$ as their underlying set. The matrices $a_*(M)$ and $b_*(N)$ are combined using pointwise sum.
\begin{defn}
For open $R$-matrices $M: X \to Y$ and $N \maps Y \to Z$ as defined above, their \define{composite} is defined by
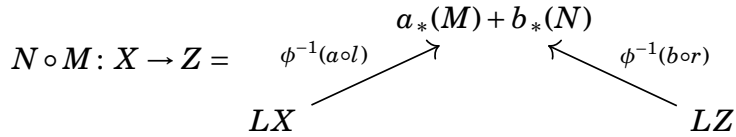
\[N \circ M \maps X \to Z = \begin{tikzcd}&a_*(M) + b_*(N)& \\
LX \ar[ur,"\phi^{-1}(a \circ l)" ] & & LZ \ar[ul,"\phi^{-1}(b \circ r)",swap] \end{tikzcd} \]
where $\phi^{-1}$ gives the unique morphism out of a discrete $R$-matrix defined by a function on its underlying set.
\end{defn}

An $R$-matrix $M \maps X \times X \to R$ can represent a graph with vertex set $X$ weighted in $R$. Similarly, an open $R$-matrix, represents an $R$-weighted graph equipped with inputs and outputs. For example, the $[0,\infty]$-matrix
\[
\begin{bmatrix}
1 &2 & .1 \\
3 & 0 & .2 \\
\infty & 1 & .2 
\end{bmatrix}
\]
on the set $\{a,b,c\}$ can be regarded as on open $[0,\infty]$-matrix with left input set $\{1,2\}$ and right input set $\{3\}$. The mappings of the cospan are given by $1 \mapsto a, 2 \mapsto b$ and $3 \mapsto c$. This can be drawn as an open weighted graph
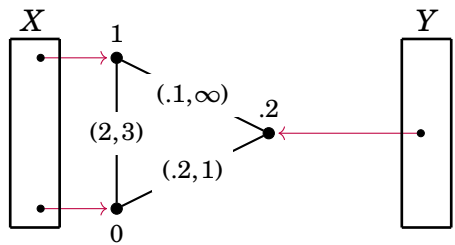
\[
\begin{tikzpicture}
	\begin{pgfonlayer}{nodelayer}
	    
		\node [style=empty] (X) at (-5.1, 2) {$X$};
		\node [style=none] (Xtr) at (-4.75, 1.75) {};
		\node [style=none] (Xbr) at (-4.75, -0.75) {};
		\node [style=none] (Xtl) at (-5.4, 1.75) {};
             \node [style=none] (Xbl) at (-5.4, -0.75) {};
	
		\node [style=inputdot] (1) at (-5, 1.5) {};
		\node [style=inputdot] (5) at (-5,-.5) {};

		\node [style=empty] (Y) at (0.1, 2) {$Y$};
		\node [style=none] (Ytr) at (.4, 1.75) {};
		\node [style=none] (Ytl) at (-.25, 1.75) {};
		\node [style=none] (Ybr) at (.4, -0.75) {};
		\node [style=none] (Ybl) at (-.25, -0.75) {};

		\node [style=inputdot] (3) at (0, 0.5) {};

		\node [style=dot,label={\footnotesize $1$}] (A) at (-4, 1.5) {};
		\node [style=dot,label=below:{\footnotesize $0$}] (A2) at (-4,-.5) {};
		\node [style=dot,label={\footnotesize $.2$}] (A3) at  (-2,.5) {};
		
	\end{pgfonlayer}
	\begin{pgfonlayer}{edgelayer}
		\draw [style=inputarrow] (1) to (A);

		\draw [style=inputarrow] (3) to (A3);
		
		\draw [style=inputarrow] (5) to (A2);
		\draw[thick] (A) -- (A2) node [midway
		, fill=white] {\footnotesize $(2,3)$};
		\draw[thick] (A) -- (A3) node [midway, fill=white] {\footnotesize $(.1,\infty)$};
		\draw[thick] (A3) -- (A2) node [midway, fill=white] {\footnotesize $(.2,1)$};
	
		\draw [style=simple] (Xtl.center) to (Xtr.center);
		\draw [style=simple] (Xtr.center) to (Xbr.center);
		\draw [style=simple] (Xbr.center) to (Xbl.center);
		\draw [style=simple] (Xbl.center) to (Xtl.center);
		\draw [style=simple] (Ytl.center) to (Ytr.center);
		\draw [style=simple] (Ytr.center) to (Ybr.center);
		\draw [style=simple] (Ybr.center) to (Ybl.center);
		\draw [style=simple] (Ybl.center) to (Ytl.center);
	\end{pgfonlayer}
\end{tikzpicture}
\]
where a tuple labeling an edge indicates the weights on that edge in both directions. Similarly, we define an open $[0,\infty]$-matrix
on $\{d,e\}$ 
\[
\begin{bmatrix}
6 & \infty\\
0& 9
\end{bmatrix}
\]
with left input set given by $\{3\}$ and right input set given by $\{4\}$. The mappings in the cospan for this open $[0,\infty]$-matrix are given by the assignments $3 \mapsto d$ and $4 \mapsto e$. This open $[0,\infty]$-matrix is drawn as
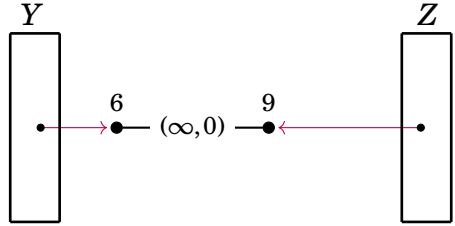
\[
\begin{tikzpicture}
	\begin{pgfonlayer}{nodelayer}
	    
		\node [style=empty] (X) at (-5.1, 2) {$Y$};
		\node [style=none] (Xtr) at (-4.75, 1.75) {};
		\node [style=none] (Xbr) at (-4.75, -0.75) {};
		\node [style=none] (Xtl) at (-5.4, 1.75) {};
             \node [style=none] (Xbl) at (-5.4, -0.75) {};
	
		\node [style=inputdot] (1) at (0, .5) {};

		\node [style=empty] (Y) at (0.1, 2) {$Z$};
		\node [style=none] (Ytr) at (.4, 1.75) {};
		\node [style=none] (Ytl) at (-.25, 1.75) {};
		\node [style=none] (Ybr) at (.4, -0.75) {};
		\node [style=none] (Ybl) at (-.25, -0.75) {};
		\node [style=inputdot] (3) at (-5, 0.5) {};
		\node [style=dot,label={\footnotesize $9$}] (A) at (-2, .5) {};
		\node [style=dot,label={\footnotesize $6$}] (A3) at  (-4,.5) {};

	\end{pgfonlayer}
	\begin{pgfonlayer}{edgelayer}
		\draw [style=inputarrow] (1) to (A);

		\draw [style=inputarrow] (3) to (A3);

		\draw[thick] (A) -- (A3) node [midway, fill=white] {\footnotesize $(\infty,0)$};

		\draw [style=simple] (Xtl.center) to (Xtr.center);
		\draw [style=simple] (Xtr.center) to (Xbr.center);
		\draw [style=simple] (Xbr.center) to (Xbl.center);
		\draw [style=simple] (Xbl.center) to (Xtl.center);
		\draw [style=simple] (Ytl.center) to (Ytr.center);
		\draw [style=simple] (Ytr.center) to (Ybr.center);
		\draw [style=simple] (Ybr.center) to (Ybl.center);
		\draw [style=simple] (Ybl.center) to (Ytl.center);
	\end{pgfonlayer}
\end{tikzpicture}
 \]
The composite of these two $[0,\infty]$-matrices is represented by
 
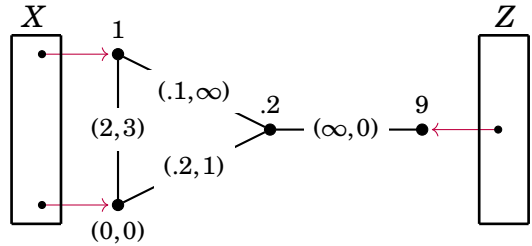
\[
\begin{tikzpicture}
	\begin{pgfonlayer}{nodelayer}
	    
		\node [style=empty] (X) at (-5.1, 2) {$X$};
		\node [style=none] (Xtr) at (-4.75, 1.75) {};
		\node [style=none] (Xbr) at (-4.75, -0.75) {};
		\node [style=none] (Xtl) at (-5.4, 1.75) {};
             \node [style=none] (Xbl) at (-5.4, -0.75) {};
	
		\node [style=inputdot] (1) at (-5, 1.5) {};
		\node [style=inputdot] (5) at (-5,-.5) {};
		\node [style=empty] (Y) at (1.1, 2) {$Z$};
		\node [style=none] (Ytr) at (1.4, 1.75) {};
		\node [style=none] (Ytl) at (.75, 1.75) {};
		\node [style=none] (Ybr) at (1.4, -0.75) {};
		\node [style=none] (Ybl) at (.75, -0.75) {};

		\node [style=inputdot] (3) at (1, 0.5) {};

		\node [style=dot,label={\footnotesize $1$}] (A) at (-4, 1.5) {};
		\node [style=dot,label=below:{\footnotesize $(0,0)$}] (A2) at (-4,-.5) {};
		\node [style=dot,label={\footnotesize $.2$}] (A3) at  (-2,.5) {};
		\node [style=dot,label={\footnotesize $9$}] (A4) at  (0,.5) {};
		
	\end{pgfonlayer}
	\begin{pgfonlayer}{edgelayer}
		\draw [style=inputarrow] (1) to (A);

		\draw [style=inputarrow] (3) to (A4);
		
		\draw [style=inputarrow] (5) to (A2);
		\draw[thick] (A) -- (A2) node [midway
		, fill=white] {\footnotesize $(2,3)$};
		\draw[thick] (A) -- (A3) node [midway, fill=white] {\footnotesize $(.1,\infty)$};
		\draw[thick] (A3) -- (A2) node [midway, fill=white] {\footnotesize $(.2,1)$};
		\draw[thick] (A3) -- (A4) node [midway, fill=white] {\footnotesize $(\infty,0)$};
	
		\draw [style=simple] (Xtl.center) to (Xtr.center);
		\draw [style=simple] (Xtr.center) to (Xbr.center);
		\draw [style=simple] (Xbr.center) to (Xbl.center);
		\draw [style=simple] (Xbl.center) to (Xtl.center);
		\draw [style=simple] (Ytl.center) to (Ytr.center);
		\draw [style=simple] (Ytr.center) to (Ybr.center);
		\draw [style=simple] (Ybr.center) to (Ybl.center);
		\draw [style=simple] (Ybl.center) to (Ytl.center);
	\end{pgfonlayer}
\end{tikzpicture}
\]
where edges are omitted if their weight is infinite in both directions. The matrix on the apex of this composite is computed by pushing each component matrix forward to the pushout of their underlying sets and adding them together i.e. 
\[
\begin{bmatrix}
1 &2 & .1 & \infty \\
3 & 0 & .2 & \infty  \\
\infty & 1 & .2 & \infty \\
\infty & \infty & \infty & \infty 
\end{bmatrix}
+
\begin{bmatrix}
\infty & \infty & \infty & \infty \\
\infty & \infty & \infty & \infty \\
\infty & \infty & 6 & \infty\\
\infty & \infty & 0 & 9
\end{bmatrix}
=
\begin{bmatrix}
1 & 2 & .1 & \infty\\
3 & 0 & .2 & \infty \\
\infty& 1 & .2 & \infty \\
\infty &\infty & 0 & 9
\end{bmatrix}
\]
The entries of this matrix represent the shortest distance between pairs of vertices.

\begin{thm}\label{open}
For a quantale $R$, there is a symmetric monoidal double category $\Open(\RMat)$ where 
\begin{itemize}
    \item objects are sets $X$,$Y$,$Z \ldots$
    \item vertical morphisms are functions $f: X \to Y$, 
    \item a horizontal morphism $M \maps X \to Y$ is an open $R$-matrix
    
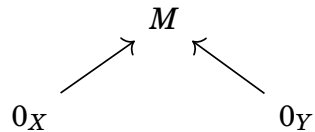
\[
    \begin{tikzcd}
     & M & \\\
    0_X \ar[ur] & & 0_Y \ar[ul]
    \end{tikzcd}
    \]
    \item vertical 2-morphisms are commutative rectangles
    
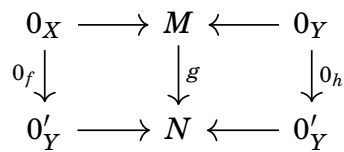
\[
    \begin{tikzcd}
    0_X\ar[d,"0_f",swap] \ar[r] & M \ar[d,"g"] & \ar[l] \ar[d,"0_h"] 0_Y \\
    0_Y' \ar[r] & N & \ar[l] 0_Y'
    \end{tikzcd}
    \]
    \item vertical composition is ordinary composition of functions,
    \item and horizontal composition is given by the composite operation defined above.
    
    
\end{itemize}
The symmetric monoidal structure is given by 
\begin{itemize}
    \item coproducts in $\Set$ on objects and vertical morphisms,
    \item pointwise coproducts on horizontal morphisms i.e. for open $R$-matrices,
    
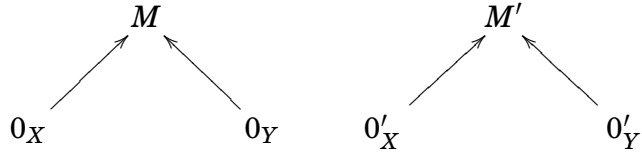
\[
    \xymatrix{ & M & & & M' & \\
	0_X \ar[ur] & & 0_Y \ar[ul] & 0_X' \ar[ur] & & 0_Y' \ar[ul] }
	\]
	their coproduct is
	
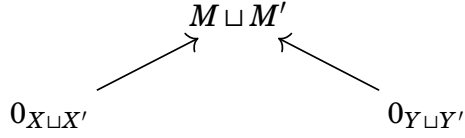
\[\begin{tikzcd}
	& M\sqcup M' & \\
	0_{X\sqcup X'} \ar[ur] & & \ar[ul] 0_{Y \sqcup Y'}
	\end{tikzcd} \]
	and pointise coproduct for two vertical 2-morphisms i.e. for vertical 2-morphisms,
	\[
	\begin{tikzcd}
    0_X\ar[d,"0_f",swap] \ar[r] & M \ar[d,"g"] & \ar[l] \ar[d,"0_h"] 0_Y \\
    0_Z \ar[r] & N & \ar[l] 0_Q
    \end{tikzcd}
    \begin{tikzcd}
    0_X'\ar[d,"0_f'",swap] \ar[r] & M' \ar[d,"g'"] & \ar[l] \ar[d,"0_h'"] 0_Y' \\
    0_Z' \ar[r] & N' & \ar[l] 0_Q'
    \end{tikzcd}
    \]
    their coproduct is 
    \[\begin{tikzcd}
    0_{X\sqcup X'}\ar[d,"0_{f\sqcup f'}",swap] \ar[r] & M \sqcup M' \ar[d,"g\sqcup g'"] & \ar[l] \ar[d,"0_{h\sqcup h'}"] 0_{Y\sqcup Y'} \\
    0_{Z\sqcup Z'} \ar[r] & N\sqcup N' & \ar[l] 0_{Q\sqcup Q'}
    \end{tikzcd} \]
\end{itemize}
\end{thm}

\begin{proof}

Theorem 3.2.3 of \cite{kenny} constructs this symmetric monoidal double category as long as
\begin{itemize}
    \item $\RMat$ has coproducts and pushouts,
    \item and $0 \maps \Set \to \RMat$ preserves pushouts and coproducts.
\end{itemize}Because $0$ is a left adjoint (Proposition \ref{L}) it preserves pushouts and coproducts when they exist so it suffices to prove the following lemma.
\end{proof}

\begin{lem}\label{matpush}
$\RMat$ has coproducts and pushouts.
\end{lem}
\begin{proof}
This is a consequence of Proposition 2.4 of \cite{wolff1974v} after noting that $\RMat$ is the category of $R$-graphs, the generating data for $R$-enriched categories. For concreteness and practicality, we offer an explicit construction of pushouts and coproducts here. 
 Let
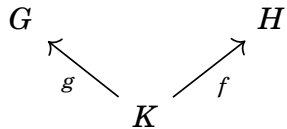
\[ \begin{tikzcd}G & & H \\
& \ar[ur,"f",swap] K \ar[ul,"g"] &\end{tikzcd}\]
be a diagram in $\RMat$ with
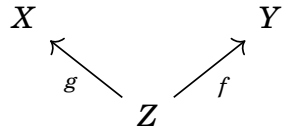
\[ \begin{tikzcd}X & & Y \\
& \ar[ur,"f",swap] Z \ar[ul,"g"] &\end{tikzcd}\]
as the underlying diagram of sets. To compute the pushout $G+_K H$ first we take the pushout of sets
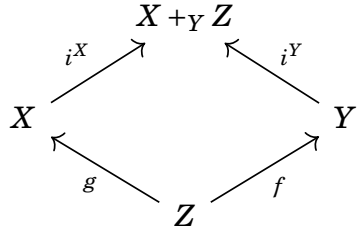
\[
\begin{tikzcd}
 & X+_Y Z & \\
X \ar[ur,"i^X"] & & Y \ar[ul,"i^Y",swap]\\
& \ar[ur,"f",swap] Z \ar[ul,"g"] &\end{tikzcd}
\]
push them forward to get matrices $i^X_*(G)$ and $i^Y_*(H)$ and join them together to get
\[G+_Y H \maps (X+_Y Z) \times (X +_Y Z) \to R = i^X_*(G) + i^Y_*(H) \]

This does indeed define a pushout in $\RMat$. Suppose we have a commutative diagram of $R$-matrices as follows:
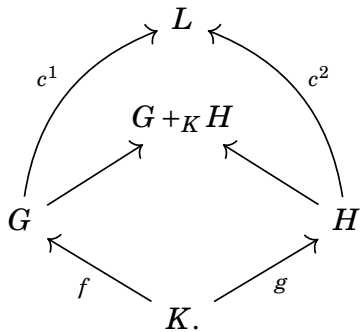
\[
\begin{tikzcd}
& L   &\\
& G+_K H   & \\
G \ar[uur,bend left,"c^1"]\ar[ur]  & & \ar[ul]\ar[uul,bend right,"c^2",swap] H \\
& K \ar[ul,"f"] \ar[ur,"g",swap]. &
\end{tikzcd}
\]
then the underlying diagram of sets induces a unique function $u$
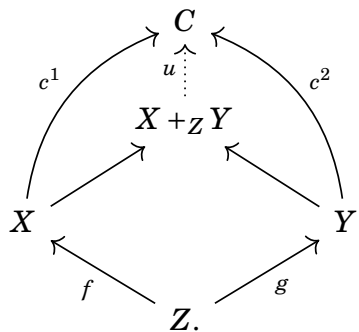
\[
\begin{tikzcd}
& C   &\\
& X+_Z Y  \ar[u,"u",dotted]  & \\
X \ar[ur]\ar[uur,bend left,"c^1"]  & & \ar[ul]\ar[uul,bend right,"c^2",swap] Y \\
& Z \ar[ul,"f"] \ar[ur,"g",swap]. &
\end{tikzcd}
\]
\noindent commuting suitable with $c^1$ and $c^2$. The map $u$ is certainly unique, it remains to show that it is well-defined i.e. it satisfies the inequality
\[u_*(G+_K H) \leq L \]
Indeed, for $(x,y) \in C \times C$,
\begin{align*}
    u_*(G+_K H)(x,y) & = \sum_{(a,b) \in (u \times u)^{-1}(x,y)} G+_K H(a,b) \\
    & = \sum_{(a,b) \in (u \times u)^{-1}(x,y)} i^X_*(G)(a,b) + i^Y_*(H) (a,b) \\
    & = \sum_{(a,b) \in (u \times u)^{-1}(x,y)} i^X_*(G)(a,b) + \sum_{(a,b) \in (u \times u)^{-1}(x,y)} i^Y_*(H) (a,b) \\
    & = u_*(i^X_*(G))(x,y) + u_*(i^Y_*(H))(x,y)
\end{align*}
However, because 
\[u_*(i^X(G)) = c^1_*(G) \text{ and } u_*(i^Y(H)) = c^2_*(G) \]
the above expression is equal to
\[c^1_*(G)(x,y) + c^2_*(H) (x,y) \]
which is less than or equal to $L(x,y)$ because each term is and $+$ is the least upper bound.

For $R$-matrices $G \maps X \times X \to R$ and $H \maps Y \times Y \to R$, their coproduct is given by the pushout

\[
\begin{tikzcd}
& G+_{\phi} H & \\
G \ar[ur]& & H\ar[ul] \\
 & \phi \ar[ul,"!_G"] \ar[ur,"!_H",swap]& 
\end{tikzcd}
\]
where $\phi$ is the unique $R$-matrix on the empty set and $!_G$ and $!_H$ are the unique morphisms into $G$ and $H$ respectively.

\end{proof}

\section{Compositionality of the Algebraic Path Problem}\label{openalg} 
In this section we show how the algebraic path problem functor $ F \maps \RMat \to \RCat$ extends to a symmetric monoidal double functor \[\Open(F) \maps \Open(\RMat) \to \Open(\RCat).\] This double functor describes how the syntax of gluing open $R$-matrices extends to a series of coherent compositionality laws for the algebraic path problem. 

For composable open $R$-matrices
\[ \begin{tikzcd}& M & & N &\\
0_X \ar[ur] & & 0_Y \ar[ul] \ar[ur] & & 0_Z \ar[ul]\end{tikzcd} \]
we apply the algebraic path problem functor $F$ to get a cospans of $R$-categories
\[ \begin{tikzcd}& F(M) & & F(N) &\\
1_X \ar[ur] & & 1_Y \ar[ul] \ar[ur] & & 1_Z \ar[ul]\end{tikzcd} \] The pushout in $\RMat$, $F(M) +_{1_Y} F(N)$, is not equal to the solution $F(M+_{0_Y} N)$. The former optimizes over only paths which start in $M$ and end in $N$. On the other hand, $F(M+_{0_Y} N)$ optimizes over paths which may zig-zag back and forth between $M$ and $N$, as many times as they like, before arriving at their destination. Therefore, to construct $F(M+_{0_Y} N)$ from its components we turn to the pushout in $\RCat$.

\begin{prop}
$\RCat$ has pushouts and coproducts.
\end{prop}

\begin{proof}
More generally, $\RCat$ has all colimits by Corollary 2.14 of \cite{wolff1974v}. These colimits are constructed via the transfinite construction of free algebras \cite{kelly1980unified}. The idea behind the transfinite construction is that colimits in a category of monoids can be constructed by first taking the colimit of their underlying objects, taking the free monoid on that colimit, and then quotienting out by the equations in your original monoids. Here we provide an explicit description in the case of $R$-categories. 
\end{proof}

\begin{prop}\label{colim}
For a diagram $D \maps C \to \RCat$, its colimit is given by the formula
\[\mathrm{colim}_{c \in C} D(c) \cong F( \mathrm{colim}_{c \in C} U(D(c)) ) \]
\end{prop}
\begin{proof}
It suffices to show that $F(\colim_{c\in C} U (D(c)))$ satisfies the universal property of $\colim_{c\in C} D(c)$. Let $\alpha \maps \Delta_d \Rightarrow D$ be a cocone from an object $d \in \RCat$ to our diagram $D$. Because $\alpha$ can be regarded as a cocone in $\RMat$, the universal property of colimits induces a unique map
\[ \colim_{c \in C} U (D(c)) \to U(d)\]
of $R$-matrices. Applying $F$ to this morphism gives a map
\[ F(\colim_{c \in C} U (D(c))) \to FU(d) = d\]
where the last equality follows either from elementary considerations or from the adjunction $F \dashv U$ being idempotent. The above map is a unique morphism satsifying the universal property for $\colim_{c\in C} D(c)$.
\end{proof}

\begin{cor}\label{comp}
For a diagram
\[
\begin{tikzcd}
M & & N \\
    & K \ar[ul] \ar[ur]& 
    \end{tikzcd}\]
    in $\RCat$, the pushout is given by
    \[M+_K N \cong F(U(M) +_{U(K)} U(N)) \]
    and the coproduct of $R$-categories is given by their coproduct in $\RMat$ i.e.
    \[M\sqcup N \cong U(M) \sqcup U(N) \]
\end{cor}
This pushout forms the horizontal composition of a double category of open $R$-categories. 
\begin{thm}There is a symmetric monoidal double category $\Open(\RCat)$ where
\begin{itemize}
\item objects are sets,
\item vertical morphisms are functions,
\item horizontal morphisms are cospans 

\[ \begin{tikzcd}
 & M & \\
   1_X \ar[ur] & & 1_Y \ar[ul]
\end{tikzcd}\]
where the apex $M$ satisfies the axioms of an $R$-category,
\item and vertical 2-morphisms are commuting rectangles
\[
\begin{tikzcd}
    1_X\ar[d,"1_f",swap] \ar[r] &  M \ar[d,"g"] & \ar[l] \ar[d,"1_h"] 1_Y \\
    1_X' \ar[r] & N & \ar[l] 1_Y'
    \end{tikzcd}
\]
\item The horizontal composition is given by pushout of open $R$-categories i.e. for open $R$-categories
\[
\begin{tikzcd}
& M & & N & \\
1_X \ar[ur] & & 1_Y \ar[ul] \ar[ur] & &\ar[ul] 1_Z
\end{tikzcd}\]
their pushout is the cospan
\[
\begin{tikzcd}
& F(U(M)+_{U(K)} U(N)) & \\
1_X \ar[ur] & & \ar[ul] 1_Y
\end{tikzcd}
\]
\end{itemize}
The symmetric monoidal structure of $\Open(\RCat)$ is given by
\begin{itemize}
    \item coproduct of sets and functions,
    \item pointwise coproduct on horizontal morphisms,
    \item and pointiwise coproduct on vertical 2-morphisms.
\end{itemize}

\end{thm}

\begin{proof}
To construct the desired symmetric monoidal double category, we apply Corollary 2.4 of \cite{structured} to the composite left adjoint \[\begin{tikzcd}\Set \ar[r,"0"] & \RMat \ar[r,"F"] & \RCat \end{tikzcd}\]
\end{proof}
\noindent So far we have the commutative diagram of functors
\[
\begin{tikzcd}
\RMat  \ar[rr,"F"] &  &\RCat \\
&\Set \ar[ur,"1",swap] \ar[ul,"0"]&
\end{tikzcd}
\]
where $1 \maps \Set \to \RCat$ is the functor which sends a set $X$ to the identity matrix $1_X$. The definition of $\Open$ is functorial with respect to this sort of diagram i.e. it induces a symmetric monoidal double functor between the relevant double categories.
\begin{thm}\label{functor}
There is a symmetric monoidal double functor
\[\bigstar \maps \Open(\RMat) \to \Open(\RCat) \]
which is 
\begin{itemize}
\item the identity on objects and vertical morphisms,

\item an open $R$-matrix
\[
M \maps X \to Y = 
\begin{tikzcd}
& M & \\
0_X \ar[ur]& & \ar[ul] 0_Y
\end{tikzcd}
\]
is sent to the solution of its algebraic path problem
\[
\bigstar(M) \maps X \to Y
\begin{tikzcd}
& FM & \\
1_X \ar[ur]& & \ar[ul] 1_Y.
\end{tikzcd}
\]

and,
\item a vertical 2-morphism of open $R$-matrices 
\[
 \alpha \maps M \Rightarrow N =	\begin{tikzcd}
    0_X\ar[d,"0_f",swap] \ar[r] & M \ar[d,"g"] & \ar[l] \ar[d,"0_h"] 0_Y \\
    0_Z \ar[r] & N & \ar[l] 0_Q
    \end{tikzcd}
\] 
is sent to the 2-morphism given by pointwise application of $F$
\[ \bigstar(\alpha) \maps M \Rightarrow N =
	\begin{tikzcd}
    1_X\ar[d,"1_f",swap] \ar[r] & FM \ar[d,"Fg"] & \ar[l] \ar[d,"1_h"] 1_Y \\
    1_X' \ar[r] & FN & \ar[l] 1_Y'
    \end{tikzcd}
\] 
\end{itemize}
\end{thm}

\begin{proof}
Theorem 4.3 of \cite{structured} proves functoriality of the ``Open" construction on squares
\[ 
\begin{tikzcd}
X\ar[r,"F_1"]  & X'\\
A \ar[r,"F_0",swap] \ar[u,"L"] & A' \ar[u,"L'",swap]
\end{tikzcd}
\]
commuting up to natural isomorphism. The result follows from applying this result to the square
\[\begin{tikzcd}\RMat \ar[r,"F"] & \RCat\\
\Set \ar[u,"0"] \ar[r,equals] & \Set \ar[u,"1",swap]
\end{tikzcd}
\]
\end{proof}\noindent The definition of symmetric monoidal double functor packages up a lot of information very succinctly. In particular, it contains coherent comparison isomorphism relating the solution of the algebraic path problem on a composite matrix to the solution on its components. For open $R$-matrices $M: X \to Y$ and $N : Y \to Z$, there is a composition comparison
\begin{equation}\label{funct}\phi_{MN} \maps \bigstar(M)\circ \bigstar(N)  \xrightarrow{\sim} \bigstar(M \circ N)\end{equation}
and monoidal comparison
\begin{equation}\label{functmon}
    \psi_{MM'} \maps \bigstar(M + M') \xrightarrow{\sim} \bigstar(M) + \bigstar(M')
\end{equation}
%
giving recipes to break solutions to the algebraic path problem into their components. In other words, the left-hand side of each comparison is computed to determine the right-hand side

Pouly and Kohlas present a similar relationship in the context of valuation algebras. \cite[\S 6.7]{pouly2012generic}. For matrices $M$ and $N$ representing weighted graphs on vertex sets $s$ and $t$ respectively, the solution to the algebraic path problem on the union of their vertex sets is given by
\[F(M) \otimes F(N) = F\left( F(M)^{\uparrow s \cup t} + F(N)^{\uparrow s \cup t} \right)\]
In this formula, $\uparrow s \cup t$ indicates that the matrix is trivially extended to the union of the vertex sets. This formula is less general than comparison (\ref{funct}): it corresponds to the special case when the legs of the open $R$-matrices are inclusions.

A typical algorithm for the algebraic path problem has spacial complexity $\Theta(n^3)$ where $n$ is the number of vertices in your weighted graph \cite{hofner2012dijkstra}. The comparisons (\ref{funct}) and (\ref{functmon}) suggest a strategy for computing the solution to the algebraic path problem which reduces this complexity. First cut your weighted graph into smaller chunks, compute the solution to the algebraic path problem on those chunks, then combine their solutions using (\ref{funct}) and \ref{functmon}). Unfortunately, this strategy will in general to take \emph{more} time to compute the solution to the algebraic path problem on a composite because the right hand side of comparison (\ref{funct}) requires three applications of the functor $F$. However, the situation improves if the open $R$-matrices are functional.

\section{Functional Open Matrices}\label{functional}
In this section we define functional open $R$-matrices, a class of open $R$-matrices for which the composition comparison
\[\phi_{MN} \maps \bigstar(M) \circ \bigstar(N) \cong \bigstar(M \circ N) \]
can be expressed in terms of matrix multiplication. The one caveat is that this expression requires that the open matrices be restricted to their inputs and outputs.
\begin{defn}
Let $M \maps X \to Y$ be the open $R$-category
\[\begin{tikzcd}
& M & \\
1_X \ar[ur,"i"] & & \ar[ul,"o",swap] 1_Y 
\end{tikzcd} \] Then the \define{blackbox} of $M$ is the matrix
\[\blacksquare (M) \maps X \times Y \to R \]
given by 
\[\blacksquare(M) ( x,y) = M(i(x),o(y)) \]
\end{defn}\noindent The $\blacksquare$ operation is extended to all of $\RCat$ but the composition is only preserved laxly.
\begin{defn}
Let $\Mat_R$ be the double category where
\begin{itemize}
    \item an object is a set $X$,$Y$,$Z$,$\ldots$
    \item a vertical morphism is a function $f\maps X \to Y$,
    \item a horizontal morphism $M\maps X \to Y$ is a matrix $M \maps X \times Y \to R$,
    \item a vertical 2-morphism from $M \maps X \to Y$ to $N \maps X' \to Y'$ is a square
    \[
    \begin{tikzcd}
    X \ar[d,"f",swap] \ar[r,"M"] & Y \ar[d,"g"] \\
    X' \ar[r,"N",swap] & Y'
    \end{tikzcd}
    \]
    such that 
    \[ \sum_{x \in f^{-1}(x'),\, y \in g^{-1}(y')} M(x,y) \leq N(x',y')  \]
    for all $x' \in X'$ and $y' \in Y'$.
    \item Vertical composition is function composition,
    \item and horizontal composition is given by matrix multiplication.
\end{itemize}
\end{defn}\noindent In this double category, the composite of matrices $M$ and $N$ is written as the juxtaposition $MN$. Blackboxing is extended to the double category of open $R$-categories.
\begin{prop}
There is a lax double functor 
\[\blacksquare \maps \Open(\RCat) \to \Mat_R \]
which  
\begin{itemize}
\item is the identity on objects,
\item sends an open $R$-category $M\maps X \to Y$ to its blackbox $\blacksquare(M)$,
\item and sends a vertical 2-cell
\[
\begin{tikzcd}
1_X \ar[r] \ar[d,"1_f",swap]  &\ar[d,"g"] M & 1_Y\ar[d,"1_h"] \ar[l] \\
1_{X'} \ar[r] & N &\ar[l] 1_{Y'}
\end{tikzcd}
\]
to the vertical 2-cell
\[
\begin{tikzcd}
X \ar[d,"f",swap]\ar[r,"\blacksquare(M)"] & Y \ar[d,"g"] \\
X' \ar[r,"\blacksquare(N)",swap] & Y' 
\end{tikzcd}
\]
\end{itemize}
\end{prop}

\begin{proof}
First observe that this lax double functor is well-defined on 2-cells. This amounts to showing that the inequality 
\begin{equation}\label{wts} \sum_{x \in f^{-1}(x'),\, y \in h^{-1}(y')} M(i(x),j(y)) \leq N(i'(x'),j'(y'))\end{equation}
holds. Because $g$ is a morphism of $R$-matrices, we have that
\begin{equation}\label{step} \sum_{a \in g^{-1}(i'(x')),\,b\in g^{-1}(j'(y'))} M(a,b) \leq N(i'(x'),j'(y')) \end{equation}Let $M(i(x),j(y))$ be a term on the left hand side of inequality (\ref{wts}).
Then by definition, $x'=f(x)$ and $y'=h(y)$ so $ a \in g^{-1}(i'(f(x))$ and $b \in g^{-1}(j'(h(y))$. However, because we started with a 2-cell in $\Open(\RCat)$, $i' \circ f = g \circ i$ and $j' \circ h = g \circ j$ so we can rewrite inequality (\ref{step}) as 
\[ \sum_{a \in g^{-1}(g \circ i(x)),\, b \in g^{-1} (g \circ j(y))} M(a,b) \leq N(i'(x'), j'(y'))\]
The term $M(i(x),j(y))$ of the left hand side of inequality (\ref{wts}) is also a term of the left hand side of inequality (\ref{step}) so we have that 
\[M(i(x),i(y)) \leq \sum_{a \in g^{-1}(g \circ i(x)),\, b \in g^{-1} (g \circ j(y))} M(a,b) \leq N(i'(x'), j'(y'))  \]
Because each term on the left hand side of (\ref{wts}) is less than the desired quantity, the join of all the terms will be as well. Therefore the lax double functor is well-defined on 2-cells. Note that $\Mat_R$ is locally posetal i.e. for every square 
 \[
    \begin{tikzcd}
    X \ar[d,"f",swap] \ar[r,"M"] & Y \ar[d,"g"] \\
    X' \ar[r,"N",swap] & Y'
    \end{tikzcd}
    \]
    there is at most one 2-cell filling it.
This property makes it so many of the axioms in the definition of lax double functor are satisfied trivially. It suffices to show that the globular composition and identity comparisons exist. The identity morphism in $\Open(\RCat)$ on a set $X$ is the cospan
\[\begin{tikzcd}
& 1_X & \\
1_X \ar[equals,ur] &  & \ar[equals,ul,,swap]1_X\end{tikzcd}  \]
The blackbox of this cospan is equal to the identity matrix on $X$, so the identity comparison is the identity.
The composition comparison
\[\blacksquare(M) \blacksquare(N) \leq \blacksquare(M \circ N) \]
follows from the chain of inequalities
\begin{align*}
   \blacksquare(M) \blacksquare(N)& = \sum_{y \in Y} \blacksquare(M)(x,y) \blacksquare(N)(y,z) \\
   &= \sum_{y \in Y} M(i(x),j(y))N(i'(y),j'(z)) \\
   &= (M +_{1_Y} N)^2 \\
   & \leq \sum_{n geq 0} (M +_{1_Y} N)^n (i(x),j'(z)) \\
   & = \blacksquare ( M \circ N) (x,z)
\end{align*}

\end{proof}
The blackboxing operation is composed with the algebraic path problem functor to get a lax symmetric monoidal double functor
\[\Open(\RMat) \xrightarrow{\bigstar} \Open(\RCat) \xrightarrow{\blacksquare} \Mat_R \]
This lax symmetric monoidal double functor gives the solution to the algebraic path problem on an open $R$-matrix \emph{restricted to its boundaries}. It is natural to ask when this mapping is strictly functorial, as this yields a very simple compositional formula for the algebraic path problem:
\[ \blacksquare( \bigstar(M \circ N)) = \blacksquare(\bigstar(M)) \blacksquare(\bigstar(N)).\]
The double functor $\blacksquare \circ \bigstar$ is strictly functorial on functional open matrices. 
\begin{defn} Let $M \maps A \times A \to R$ be an $R$-matrix. An element $a \in X$ is a \define{source} if for every $b \in X$, $M(b,a)=0$ and a \define{sink} if $M(a,b)=0$. A \define{functional open $R$-matrix} is an open $R$-matrix 
\[
\begin{tikzcd}
& M & \\
0_X \ar[ur,"l"] & & \ar[ul,"r",swap] 0_Y
\end{tikzcd}
\]
such that for every $x \in X$, $l(x)$ is a source and for every $y \in Y$, $r(y)$ is a sink.
\end{defn}

 

\noindent Because the pushout of functional open $R$-matrices is also functional, we can form the following sub-double category.
\begin{defn}
Let $\Open(\RMat)_{fxn}$ be the full sub-symmetric monoidal double category generated by the open $R$-matrices which are functional.
\end{defn}

\begin{thm}\label{strict}
The composite $\blacksquare \circ \bigstar$ restricts to a strict double functor
\[\blacksquare \circ \bigstar_{fxn} \maps \Open(\RMat)_{fxn} \to \Mat_R \]
\end{thm}\noindent The proof of this theorem relies on a lemma which resembles the the binomial expansion of $(a+b)^n$ in a ring where $ba=0$.
\begin{lem}\label{binomial}
for functional open $R$-matrices $M \maps X \to Y$ and $N \maps Y \to Z$ we have that
\[\blacksquare (M +_{1_Y} N)^n = \sum_{i+j=n} \blacksquare(M^i)\blacksquare(N^j)\]
\end{lem}
\begin{proof}
The entries of the left hand side are expanded as 
\[\blacksquare( (M +_{1_Y} N)^n)(a_0,a_{n}) = \sum_{a_1,a_2,\ldots,a_{n-1}} (M +_{1_Y} N)(a_0,a_1) (M +_{1_Y} N)(a_1,a_2) \ldots (M +_{1_Y} N)(a_{n-1},a_n)  \]
where the $a_i$ are equivalence classes in $RM +_Y RN$. For a particular term of this sum, let $1 \leq k \leq n$ be the first natural number such that $a_k$ contains an element of $RN$. Because $M$ and $N$ are functional, for $k \leq i \leq n$ the equivalence classes $a_i$ must also contain an element of $RN$ if our term is nonzero. Therefore for a fixed $k$ the contribution to the above sum is given by
\[
\sum M(a_0,a_1) \ldots M(a_{k-1},a_k) N(a_k,a_{k+1}) \ldots N(a_{n-1},a_n)
\]
which simplifies to 
\[\blacksquare(M^k)\blacksquare( N^{n-k}) (a_0,a_n). \]
Because $k$ can occur in any entry we have that
\begin{align*}
    \blacksquare((M+_{1_Y} N)^n) & = \sum_{k \leq n} \blacksquare(M^k) \blacksquare(N^{n-k}) \\
    &= \sum_{i+j=n} \blacksquare(M^i) \blacksquare(N^j)
\end{align*}
\end{proof}\noindent \textbf{Proof of Theorem \ref{strict}:} It suffices to prove that for functional open matrices \[\begin{tikzcd}0_X \ar[r] & M & \ar[l] 0_Y \end{tikzcd}\] and \[\begin{tikzcd} 0_Y \ar[r] & N & \ar[l] 0_Z \end{tikzcd}\] the equation
\[\blacksquare( \bigstar(M \circ N)) = \blacksquare(\bigstar(M)) \blacksquare(\bigstar(N)) \]
holds. Consider the left-hand side:
\begin{align*}
   \blacksquare \bigstar M \circ N & = \blacksquare \sum_{n \geq 0} (M \circ N)^n \\
   &= \sum_{n \geq 0} \blacksquare (M \circ N)^n \\
   &=\sum_{n \geq 0} \sum_{i+j=n} \blacksquare(M^i)\blacksquare(N^j) \\
\end{align*}
on the other hand, 
\begin{align*}
    \blacksquare(\bigstar(M))\blacksquare(\bigstar(N)) &= \sum_{i \geq 0} \blacksquare(M^i) \sum_{j \geq 0} \blacksquare(N^j) \\
    &= \sum_{i,j \geq 0} \blacksquare(M^i) \blacksquare(M^j)
\end{align*}
Both sums contain the term $\blacksquare(M^i)\blacksquare(N^j)$ for every value of $i$ and $j$, but the left hand side may contain repeated terms. However, because addition is idempotent, repeated terms don't contribute to the sum and the two sides are the same.

\hfill $\square$

\section{Conclusion}
The functoriality of Theorem \ref{strict} might not be surprising. It says that if your open matrices are joined together directionally along bottlenecks, then the computation of the algebraic path problem can be reduced to a computation on components. This strategy has already proven sucessful. In \cite{sairam1992divide}, Sairam, Tamassia, and Vitter show how choosing \emph{one way separators} as cuts in a graph, allow for an efficient divide and conquer parallel algorithm for computing shortest paths. In \cite{rathke2014compositional} Rathke, Sobocinksi, and Stephens show how the reachability problem on a 1-safe Petri net can be computed more efficiently by cutting it up into more manageable pieces. Theorem \ref{functor} provides a framework for compositional formulas of this type. In future work we plan on extending the construction of this theorem to many other sorts of discrete event dynamic systems.

Lemma \ref{binomial} also holds independent computational interest. The equation given there gives a novel compositional formula for computing the solution to the algebraic path problem. The author has implemented this formula for the special case of Markov processes \cite{compmarkov}. We hope that this is the start of a more extensive library, made faster and more reliable by the mathematics developed in this paper.

\section{Acknowledgements}
I would like to thank Mike Shulman, John Baez, Christian Williams, Joe Moeller, Rany Tith, Sarah Rovner-Frydman, Zans Mihejez, Oscar Hernandez, Alex Pokorny and Todd Trimble for their helpful comments and contributions. I would also like to thank everyone in my life who supported me during this time, in particular Allison Lucas. Your work contributed to this paper as well. This work was produced on Tongva land.
\bibliographystyle{alpha}
\bibliography{petri}
\end{document}